\documentclass[11pt]{article}
\usepackage[letterpaper,margin=0.88in]{geometry}
\usepackage{amsmath,amssymb,booktabs,tikz}
\usetikzlibrary{arrows.meta}
\usepackage[hidelinks]{hyperref}
\definecolor{tre}{HTML}{214D78}
\definecolor{boxa}{HTML}{A34622}
\definecolor{boxb}{HTML}{147566}
\newcommand{\F}{\mathbb F_2}

\usepackage{amsthm}
\theoremstyle{plain}
\newtheorem{theorem}{Theorem}

\newtheorem{conjecture}{Conjecture}

\usepackage[normalem]{ulem}
\theoremstyle{remark}
\newtheorem{rem}{Remark}

\newcommand{\gauss}[3]{%
  \draw[black!70] (0,0) circle (1.36);
  \draw[-{Stealth[length=3pt]},black!75] (101:1.36) arc (101:93:1.36);
  \draw[black,line width=1pt] (90:1.29)--(90:1.43);
  \foreach \entry [count=\k from 1] in {#1}{
    \pgfmathtruncatemacro{\lab}{abs(\entry)}
    \pgfmathsetmacro{\ang}{90-360*(\k-.5)/#2}
    \ifnum\entry>0 \coordinate (tail\lab) at (\ang:1.36);
    \else \coordinate (head\lab) at (\ang:1.36);\fi
    \node[font=\scriptsize] at (\ang:1.63) {\(\ifnum\entry>0+\fi\entry\)};
  }
  \foreach \lab in {#3}{
    \ifnum\lab<4 \def\ink{tre}\else\ifnum\lab<7 \def\ink{boxa}\else\def\ink{boxb}\fi\fi
    \draw[\ink,line width=.9pt,-{Stealth[length=4.2pt,width=3.8pt]}]
       (tail\lab)--(head\lab);
  }}
\title{A counterexample to Manturov's connected-sum conjecture}
\author{V. Chernov, R. Sadykov}
\date{}

\begin{document}
\maketitle
\begin{abstract}
Manturov conjectured that his parity projection of a connected sum of two virtual knots is isotopic to the connected sum of the projections. We construct a counterexample in which a specified connected sum of two virtual knots projects to the unknot, while the connected sum of the projections is the right-handed trefoil. The counterexample is given by an explicit Gauss word with nine classical crossings, suggested by GPT-6 Astra. 
We use Turaev's primitive based matrices to show that the diagrams of the summands and of their sum all have minimal supporting genus, as Manturov’s parity projection construction requires. 
\end{abstract}
\medskip
\noindent\textit{2020 Mathematics Subject Classification.} Primary 57K12; Secondary 57K10.

\noindent\textit{Keywords.} Virtual knots, parity, connected sum, Gauss words.

\section{Introduction}

An oriented {\em virtual knot diagram} is a generic immersion of an oriented circle into the plane
whose double points are either classical crossings, with an overpassing
and an underpassing branch, or virtual crossings, marked by small
circles. Two diagrams represent the same {\em virtual knot} if they are
related by isotopy, the usual Reidemeister moves at classical
crossings, and detour moves. A detour move replaces an arc meeting only
virtual crossings by another arc with the same endpoints; all new
intersections are declared virtual \cite{Kauffman}.
To form a connected sum of two virtual knots, choose a point away from all the crossings in
each of two virtual knot diagrams, cut the diagrams at these points, and join the terminal endpoint of each arc to the initial endpoint of the other. All new intersections are
declared virtual. The resulting virtual knot can depend on the
diagrams and the cuts \cite{KishinoSatoh}, so these choices are part of the construction.

 A virtual knot is {\em classical} if it
admits a diagram with no virtual crossings. Goussarov, Polyak, and Viro\cite{GPV} showed that two classical diagrams that are equivalent as virtual knots are equivalent as classical knots, so there is an embedding of classical knot types into virtual knot types. The same conclusion follows from Kuperberg's uniqueness theorem  \cite{Kuperberg}.  
Manturov constructed a projection $g$ from virtual knots to classical knots which fixes every classical knot \cite[Theorem 4 and Sec. 4.1]{M}. Chernov and Sadykov constructed an analogous projection from virtual Legendrian knots in $ST^*F$ to Legendrian knots in $ST^*S^2$ \cite{ChernovSadykov}.
Manturov conjectured that $g$ preserves connected sums.

\begin{conjecture} [Manturov {\cite[Sec.~4.1, Conjecture~1]{M}}]For every specified connected sum of virtual knots $K_1$ and $K_2$,
\[
g(K_1\#K_2)=g(K_1)\#g(K_2),
\]
where the sum on the right is the ordinary connected sum of classical
knots, and equality means that the two classical knots are isotopic.
\end{conjecture}

The main result of the paper is the existence of a counterexample. 

\begin{theorem}\label{th:main}
    There are virtual knots $K_1$ and $K_2$, and a connected sum $K=K_1\#K_2$ such that $g(K_1)\#g(K_2)$ is the right-handed trefoil, while $g(K)$ is the unknot. In particular, Manturov's conjecture is false. 
\end{theorem}

We construct the counterexample starting with the standard diagram of the right-handed trefoil. The second summand $K_2$ is a non-classical virtual knot, and the first summand $K_1$ is a connected sum of $K_2$ with the trefoil. To obtain the diagram of $K_1$, we replace an arc of the trefoil with the complement of an arc in a diagram of $K_2$. We form the connected sum $K$ by replacing a different arc of the trefoil in $K_1$ with the complement of an arc in the diagram of $K_2$. 

In the diagram of $K_1$ in Figure~\ref{fig:virtual-summands}, the three trefoil crossings $1, 2$ and $3$ are even with respect to Manturov's parity group (see \S\ref{sec:manturov}), while the crossings $4, 5$ and $6$ are odd. 
The parity step makes the odd crossings $4, 5$ and $6$ virtual, and detour moves then remove the inserted copy of $K_2$, leaving the standard trefoil diagram. Hence $g(K_1)$ is the right-handed trefoil. 
In the diagrams of $K_2$ and $K$, every classical crossing is odd. For $K$ this includes the trefoil crossings $1, 2$ and $3$, which are even in $K_1$: the diagram of $K$ contains a second copy of $K_2$ inserted into a different arc of the trefoil, and this changes the parity of these crossings. The parity step makes all crossings virtual, and then detour moves turn each diagram into a simple closed curve. Hence $g(K_2)$ and $g(K)$ are both the unknot. 
The projection may be computed only on diagrams of minimal supporting genus. To justify applying it to diagrams of $K_1, K_2$ and $K$, we use Turaev's primitive based matrices \cite{T} to show that their minimal supporting genera are $2,2$ and $4$, respectively.

The Gauss word of the counterexample was suggested by GPT-6 Astra (OpenAI). GPT-6 Astra
and Claude Opus~5.5 (Anthropic) were also used substantially in drafting
and editing the paper. The authors  take
full responsibility for the content.

\section{Manturov's projection}\label{sec:manturov}

Let $D$ be a virtual knot diagram. Regard its classical crossings as
vertices and the arcs between consecutive classical crossings as edges.
Retain the cyclic order of the four branches at each vertex. The graph
with these cyclic orders is the \emph{ribbon graph} of $D$. Replace its
vertices by disks and its edges by bands, attached in the prescribed
cyclic order and compatibly with orientation. At virtual crossings the
bands remain disjoint. Capping the boundary components with disks gives
a closed oriented surface $F_D$ on which $D$ has only classical
crossings. A diagram without classical crossings is thickened to an
annulus and capped to a sphere. The \emph{minimal supporting genus}
of a virtual knot is the least genus of $F_D$ over all its diagrams. By Kuperberg's theorem~\cite{Kuperberg}, a representative of a virtual knot in a thickened surface of minimal genus is unique up to diffeomorphism and isotopy of the thickened surface. 

Write $s=[D]\in H_1(F_D;\F)$ for the homology class of the whole curve. Following Ilyutko, Manturov and Nikonov~\cite{IMN}, 
Manturov identifies the \emph{parity group} of $D$ with
$G(D)=H_1(F_D;\F)/\langle s\rangle$ \cite[Theorem~7]{M}. At a classical crossing $v$, smoothing according to
the orientation splits the curve into two closed curves, called its
{\em halves.} If one half has homology class $h_v$, the other has class
$s+h_v$. Their images in $G(D)$ therefore agree. The crossing is
\emph{even} if this common image is zero, and \emph{odd} otherwise.
Thus $v$ is even precisely when $h_v\in\{0,s\}$, or equivalently, when
one of its halves is null-homologous in $F_D$. A {\em parity step\/} replaces
all odd crossings by virtual crossings.

Let $K$ be a virtual knot. Manturov's {\em projection map} $g$ \cite[Sec. 4.1]{M} is obtained by iterating parity steps. Choose a diagram of $K$ of minimal supporting genus, and apply a parity step. Repeat this with the resulting knot until the resulting knot is classical,  recomputing the parity
group at each step. The result does not depend on the diagrams chosen, so the procedure defines a map $g$ from virtual knots to classical knots which fixes classical knots \cite[Theorem~4, Theorem 5 and Sec.~4.1]{M}.

\section{The two summands}\label{sec:summands}

In this section we construct two virtual knots $K_1$ and $K_2$ and a connected sum $K=K_1\#K_2$  for which Manturov's conjecture fails. 

We describe the knot diagrams by Gauss words. Label the classical crossings
of an oriented virtual knot diagram by distinct positive integers.
Starting at a marked point and following the orientation of the knot,
record $+i$ when passing over at crossing $i$, and $-i$ when passing
under. Each label occurs exactly twice, once with each sign, and the word is read cyclically: its last letter is followed by its first. {\em The signs in the word indicate overpasses and underpasses, not the signs of the crossings.}
A Gauss word together with the signs of its crossings determines a virtual knot diagram up to detour moves \cite{Kauffman, GPV}.  In our examples, all classical crossings are positive.  For a word $X$ we write $\operatorname{Cl}(X)$ for the virtual knot represented by $X$ read cyclically. 

Consider the virtual knots $K_1$ and $K_2$ represented by the words
\[
\begin{aligned}
 W_1&=(+1,-2,+3,\; -5,-4,-6,+5,+6,+4,\; -1,+2,-3),\\
 W_2&=(-8,-7,-9,+8,+9,+7).
\end{aligned}
\]

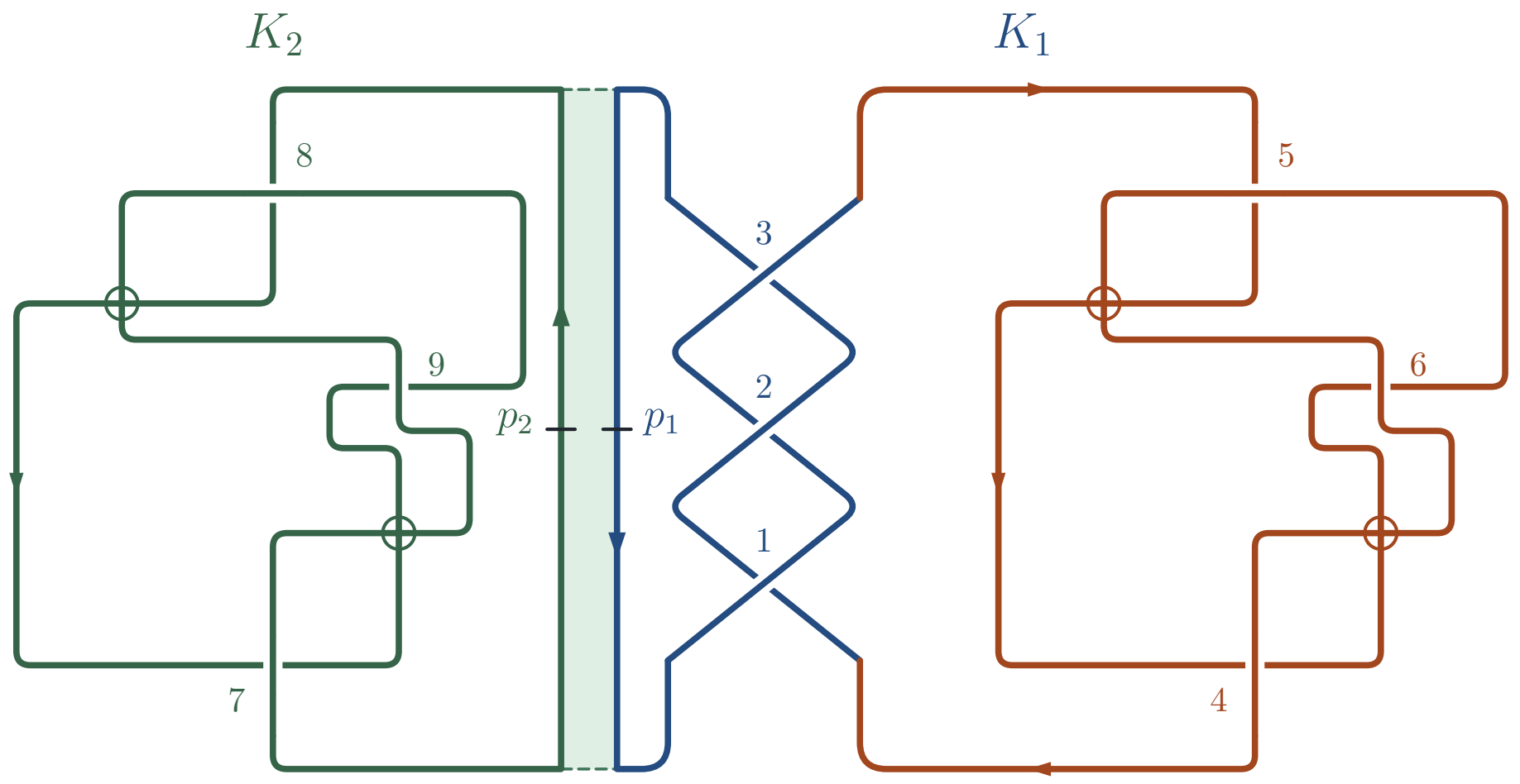
\begin{figure}[htbp]
\centering
\resizebox{0.7\linewidth}{!}{%
\begin{tikzpicture}[x=1bp,y=1bp,line cap=round,line join=round]
\definecolor{m244C81}{HTML}{244C81}
\definecolor{mA2461E}{HTML}{A2461E}
\definecolor{m366448}{HTML}{366448}
\definecolor{m20262E}{HTML}{20262E}
\definecolor{m52616A}{HTML}{52616A}
\definecolor{mD6EBDC}{HTML}{D6EBDC}
\definecolor{m40765A}{HTML}{40765A}
\definecolor{mFFFFFF}{HTML}{FFFFFF}
\fill[mD6EBDC!78] (301,53.8) rectangle (317,248.2);
\draw[m40765A,densely dashed,line width=.85bp] (301,248.2)--(317,248.2);
\draw[m40765A,densely dashed,line width=.85bp] (301,53.8)--(317,53.8);
\draw[m244C81,line width=1.80bp] (331.55,84.85)--(382.24,125.57).. controls (385.05,127.82) and (385.05,130.08) ..(382.24,132.33)--(335.76,169.67).. controls (332.95,171.92) and (332.95,174.18) ..(335.76,176.43)--(386.45,217.15);
\draw[m244C81,line width=1.80bp] (386.45,84.85)--(335.76,125.57).. controls (332.95,127.82) and (332.95,130.08) ..(335.76,132.33)--(382.24,169.67).. controls (385.05,171.92) and (385.05,174.18) ..(382.24,176.43)--(331.55,217.15);
\draw[mA2461E,line width=1.80bp] (386.45,217.15)--(386.45,240.55).. controls (386.45,245.65) and (389.00,248.20) ..(394.10,248.20)--(495.43,248.20).. controls (498.08,248.20) and (499.40,246.88) ..(499.40,244.23)--(499.40,238.75);
\draw[mA2461E,line width=1.80bp] (499.40,63.25)--(499.40,57.77).. controls (499.40,55.12) and (498.08,53.80) ..(495.43,53.80)--(394.10,53.80).. controls (389.00,53.80) and (386.45,56.35) ..(386.45,61.45)--(386.45,84.85);
\draw[mA2461E,line width=1.80bp] (499.40,238.75)--(499.40,191.05).. controls (499.40,188.35) and (498.05,187.00) ..(495.35,187.00)--(430.10,187.00).. controls (427.40,187.00) and (426.05,185.65) ..(426.05,182.95)--(426.05,87.55).. controls (426.05,84.85) and (427.40,83.50) ..(430.10,83.50)--(531.35,83.50).. controls (534.05,83.50) and (535.40,84.85) ..(535.40,87.55)--(535.40,141.55).. controls (535.40,144.25) and (534.05,145.60) ..(531.35,145.60)--(519.65,145.60).. controls (516.95,145.60) and (515.60,146.95) ..(515.60,149.65)--(515.60,159.10).. controls (515.60,161.80) and (516.95,163.15) ..(519.65,163.15)--(566.90,163.15).. controls (569.60,163.15) and (570.95,164.50) ..(570.95,167.20)--(570.95,214.45).. controls (570.95,217.15) and (569.60,218.50) ..(566.90,218.50)--(460.25,218.50).. controls (457.55,218.50) and (456.20,217.15) ..(456.20,214.45)--(456.20,180.70).. controls (456.20,178.00) and (457.55,176.65) ..(460.25,176.65)--(531.35,176.65).. controls (534.05,176.65) and (535.40,175.30) ..(535.40,172.60)--(535.40,154.60).. controls (535.40,151.90) and (536.75,150.55) ..(539.45,150.55)--(551.60,150.55).. controls (554.30,150.55) and (555.65,149.20) ..(555.65,146.50)--(555.65,125.35).. controls (555.65,122.65) and (554.30,121.30) ..(551.60,121.30)--(503.45,121.30).. controls (500.75,121.30) and (499.40,119.95) ..(499.40,117.25)--(499.40,63.25);
\draw[mFFFFFF,line width=5.50bp] (505.20,218.50)--(493.60,218.50);
\draw[mA2461E,line width=1.80bp] (507.70,218.50)--(491.10,218.50);
\node[inner sep=0pt,anchor=base,text=mA2461E,font={\fontsize{10bp}{12bp}\selectfont}] at (508.40,226.15) {5};
\draw[mFFFFFF,line width=5.50bp] (499.40,89.30)--(499.40,77.70);
\draw[mA2461E,line width=1.80bp] (499.40,91.80)--(499.40,75.20);
\node[inner sep=0pt,anchor=base,text=mA2461E,font={\fontsize{10bp}{12bp}\selectfont}] at (489.05,70.45) {4};
\draw[mFFFFFF,line width=5.50bp] (535.40,168.95)--(535.40,157.35);
\draw[mA2461E,line width=1.80bp] (535.40,171.45)--(535.40,154.85);
\node[inner sep=0pt,anchor=base,text=mA2461E,font={\fontsize{10bp}{12bp}\selectfont}] at (546.20,166.30) {6};
\draw[mA2461E,line width=.95bp] (456.20,187.00) circle[radius=4.6bp];
\draw[mA2461E,line width=.95bp] (535.40,121.30) circle[radius=4.6bp];
\fill[mA2461E,line width=1.80bp] (426.05,131.95)--(428.10,138.55)--(424.00,138.55)--cycle;
\fill[mA2461E,line width=1.80bp] (441.05,248.20)--(434.45,250.25)--(434.45,246.15)--cycle;
\fill[mA2461E,line width=1.80bp] (434.45,53.80)--(441.05,51.75)--(441.05,55.85)--cycle;
\draw[m366448,line width=1.80bp] (218.60,238.75)--(218.60,191.05).. controls (218.60,188.35) and (217.25,187.00) ..(214.55,187.00)--(149.30,187.00).. controls (146.60,187.00) and (145.25,185.65) ..(145.25,182.95)--(145.25,87.55).. controls (145.25,84.85) and (146.60,83.50) ..(149.30,83.50)--(250.55,83.50).. controls (253.25,83.50) and (254.60,84.85) ..(254.60,87.55)--(254.60,141.55).. controls (254.60,144.25) and (253.25,145.60) ..(250.55,145.60)--(238.85,145.60).. controls (236.15,145.60) and (234.80,146.95) ..(234.80,149.65)--(234.80,159.10).. controls (234.80,161.80) and (236.15,163.15) ..(238.85,163.15)--(286.10,163.15).. controls (288.80,163.15) and (290.15,164.50) ..(290.15,167.20)--(290.15,214.45).. controls (290.15,217.15) and (288.80,218.50) ..(286.10,218.50)--(179.45,218.50).. controls (176.75,218.50) and (175.40,217.15) ..(175.40,214.45)--(175.40,180.70).. controls (175.40,178.00) and (176.75,176.65) ..(179.45,176.65)--(250.55,176.65).. controls (253.25,176.65) and (254.60,175.30) ..(254.60,172.60)--(254.60,154.60).. controls (254.60,151.90) and (255.95,150.55) ..(258.65,150.55)--(270.80,150.55).. controls (273.50,150.55) and (274.85,149.20) ..(274.85,146.50)--(274.85,125.35).. controls (274.85,122.65) and (273.50,121.30) ..(270.80,121.30)--(222.65,121.30).. controls (219.95,121.30) and (218.60,119.95) ..(218.60,117.25)--(218.60,63.25);
\draw[mFFFFFF,line width=5.50bp] (224.40,218.50)--(212.80,218.50);
\draw[m366448,line width=1.80bp] (226.90,218.50)--(210.30,218.50);
\node[inner sep=0pt,anchor=base,text=m366448,font={\fontsize{10bp}{12bp}\selectfont}] at (227.60,226.15) {8};
\draw[mFFFFFF,line width=5.50bp] (218.60,89.30)--(218.60,77.70);
\draw[m366448,line width=1.80bp] (218.60,91.80)--(218.60,75.20);
\node[inner sep=0pt,anchor=base,text=m366448,font={\fontsize{10bp}{12bp}\selectfont}] at (208.25,70.45) {7};
\draw[mFFFFFF,line width=5.50bp] (254.60,168.95)--(254.60,157.35);
\draw[m366448,line width=1.80bp] (254.60,171.45)--(254.60,154.85);
\node[inner sep=0pt,anchor=base,text=m366448,font={\fontsize{10bp}{12bp}\selectfont}] at (265.40,166.30) {9};
\draw[m366448,line width=.95bp] (175.40,187.00) circle[radius=4.6bp];
\draw[m366448,line width=.95bp] (254.60,121.30) circle[radius=4.6bp];
\fill[m366448,line width=1.80bp] (145.25,131.95)--(147.30,138.55)--(143.20,138.55)--cycle;
\draw[mFFFFFF,line width=5.50bp] (354.48,103.27)--(363.52,110.53);
\draw[m244C81,line width=1.80bp] (352.53,101.70)--(365.47,112.10);
\node[inner sep=0pt,anchor=base,text=m244C81,font={\fontsize{10bp}{12bp}\selectfont}] at (359.00,115.90) {1};
\draw[mFFFFFF,line width=5.50bp] (354.48,147.37)--(363.52,154.63);
\draw[m244C81,line width=1.80bp] (352.53,145.80)--(365.47,156.20);
\node[inner sep=0pt,anchor=base,text=m244C81,font={\fontsize{10bp}{12bp}\selectfont}] at (359.00,160.00) {2};
\draw[mFFFFFF,line width=5.50bp] (354.48,191.47)--(363.52,198.73);
\draw[m244C81,line width=1.80bp] (352.53,189.90)--(365.47,200.30);
\node[inner sep=0pt,anchor=base,text=m244C81,font={\fontsize{10bp}{12bp}\selectfont}] at (359.00,204.10) {3};
\draw[m366448,line width=1.80bp] (218.60,238.75)--(218.60,244.23).. controls (218.60,246.88) and (219.92,248.20) ..(222.57,248.20)--(301,248.2)--(301,53.8)--(222.57,53.80).. controls (219.92,53.80) and (218.60,55.12) ..(218.60,57.77)--(218.60,63.25);
\draw[m244C81,line width=1.80bp] (331.55,217.15)--(331.55,240.55).. controls (331.55,245.65) and (329.00,248.20) ..(323.90,248.20)--(317,248.2)--(317,53.8)--(323.90,53.80).. controls (329.00,53.80) and (331.55,56.35) ..(331.55,61.45)--(331.55,84.85);
\fill[m366448] (301.00,188.00)--(298.50,180.50)--(303.50,180.50)--cycle;
\fill[m244C81] (317.00,114.00)--(319.50,121.50)--(314.50,121.50)--cycle;
\draw[m20262E,line width=1bp] (297,151)--(305,151);
\draw[m20262E,line width=1bp] (313,151)--(321,151);
\node[inner sep=0pt,text=m366448,font={\fontsize{12bp}{14bp}\selectfont}] at (288,153) {$p_2$};
\node[inner sep=0pt,text=m244C81,font={\fontsize{12bp}{14bp}\selectfont}] at (330,153) {$p_1$};
\node[inner sep=0pt,text=m366448,font={\fontsize{14bp}{16bp}\selectfont}] at (219,264) {$K_2$};
\node[inner sep=0pt,text=m244C81,font={\fontsize{14bp}{16bp}\selectfont}] at (433,264) {$K_1$};
\end{tikzpicture}}
\caption{The marked summands $K_2$ (left) and $K_1$ (right). The two vertical sides of the shaded rectangle contain $p_2$ and $p_1$. Replacing these sides by the dashed horizontal sides gives Figure~\ref{fig:insertions}. Small circles mark virtual crossings. All classical crossings are positive.}
\label{fig:virtual-summands}
\end{figure}

Diagrams realizing these words are shown in Figure~\ref{fig:virtual-summands}.

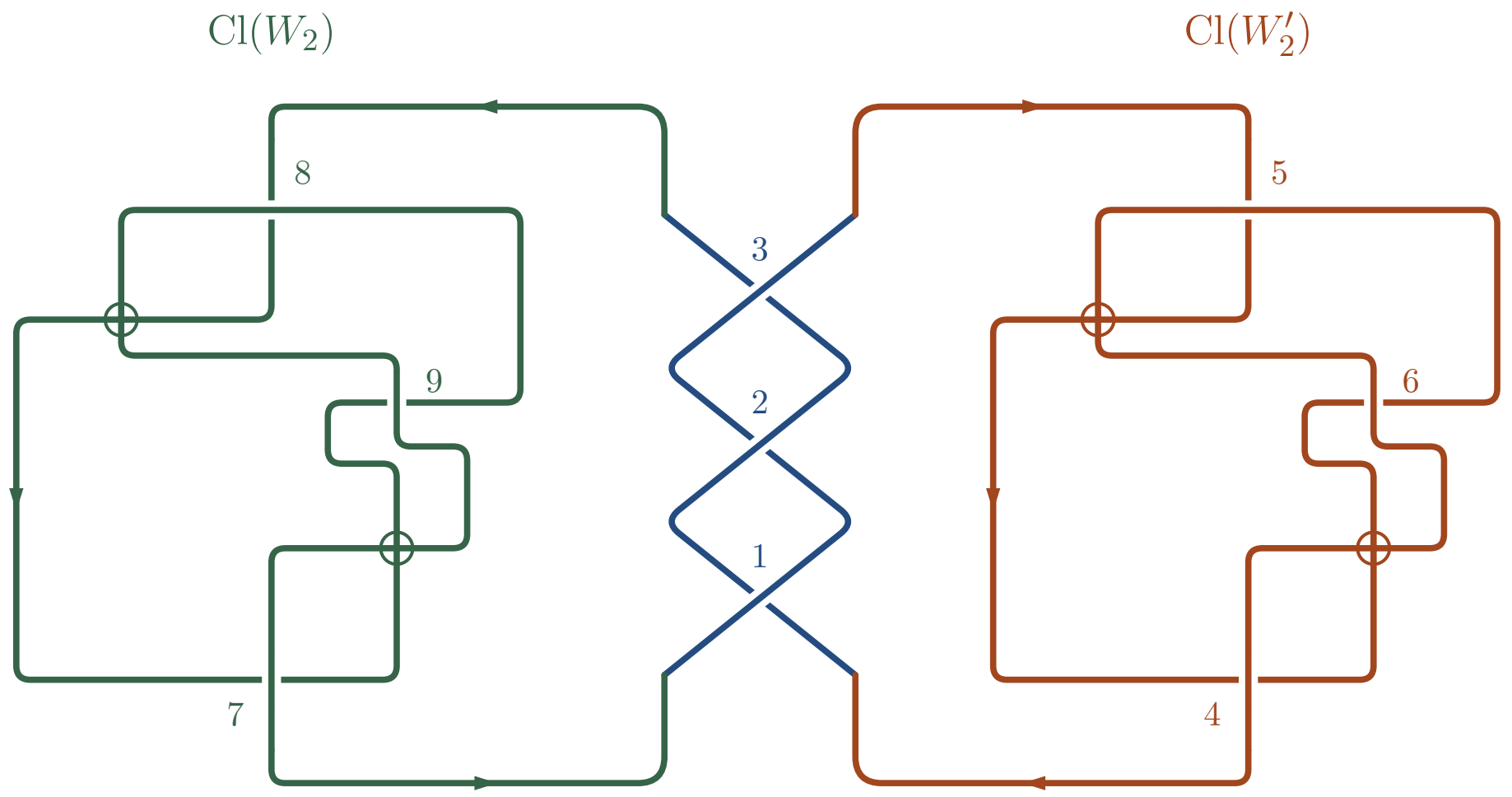
\begin{figure}[htbp]
\centering
\resizebox{0.7\linewidth}{!}{%
\begin{tikzpicture}[x=1bp,y=1bp,line cap=round,line join=round]
\definecolor{m244C81}{HTML}{244C81}
\definecolor{mA2461E}{HTML}{A2461E}
\definecolor{m366448}{HTML}{366448}
\definecolor{m20262E}{HTML}{20262E}
\definecolor{m52616A}{HTML}{52616A}
\definecolor{mD6EBDC}{HTML}{D6EBDC}
\definecolor{m40765A}{HTML}{40765A}
\definecolor{mFFFFFF}{HTML}{FFFFFF}
\draw[m244C81,line width=1.80bp] (331.55,84.85)--(382.24,125.57).. controls (385.05,127.82) and (385.05,130.08) ..(382.24,132.33)--(335.76,169.67).. controls (332.95,171.92) and (332.95,174.18) ..(335.76,176.43)--(386.45,217.15);
\draw[m244C81,line width=1.80bp] (386.45,84.85)--(335.76,125.57).. controls (332.95,127.82) and (332.95,130.08) ..(335.76,132.33)--(382.24,169.67).. controls (385.05,171.92) and (385.05,174.18) ..(382.24,176.43)--(331.55,217.15);
\draw[mA2461E,line width=1.80bp] (386.45,217.15)--(386.45,240.55).. controls (386.45,245.65) and (389.00,248.20) ..(394.10,248.20)--(495.43,248.20).. controls (498.08,248.20) and (499.40,246.88) ..(499.40,244.23)--(499.40,238.75);
\draw[mA2461E,line width=1.80bp] (499.40,63.25)--(499.40,57.77).. controls (499.40,55.12) and (498.08,53.80) ..(495.43,53.80)--(394.10,53.80).. controls (389.00,53.80) and (386.45,56.35) ..(386.45,61.45)--(386.45,84.85);
\draw[mA2461E,line width=1.80bp] (499.40,238.75)--(499.40,191.05).. controls (499.40,188.35) and (498.05,187.00) ..(495.35,187.00)--(430.10,187.00).. controls (427.40,187.00) and (426.05,185.65) ..(426.05,182.95)--(426.05,87.55).. controls (426.05,84.85) and (427.40,83.50) ..(430.10,83.50)--(531.35,83.50).. controls (534.05,83.50) and (535.40,84.85) ..(535.40,87.55)--(535.40,141.55).. controls (535.40,144.25) and (534.05,145.60) ..(531.35,145.60)--(519.65,145.60).. controls (516.95,145.60) and (515.60,146.95) ..(515.60,149.65)--(515.60,159.10).. controls (515.60,161.80) and (516.95,163.15) ..(519.65,163.15)--(566.90,163.15).. controls (569.60,163.15) and (570.95,164.50) ..(570.95,167.20)--(570.95,214.45).. controls (570.95,217.15) and (569.60,218.50) ..(566.90,218.50)--(460.25,218.50).. controls (457.55,218.50) and (456.20,217.15) ..(456.20,214.45)--(456.20,180.70).. controls (456.20,178.00) and (457.55,176.65) ..(460.25,176.65)--(531.35,176.65).. controls (534.05,176.65) and (535.40,175.30) ..(535.40,172.60)--(535.40,154.60).. controls (535.40,151.90) and (536.75,150.55) ..(539.45,150.55)--(551.60,150.55).. controls (554.30,150.55) and (555.65,149.20) ..(555.65,146.50)--(555.65,125.35).. controls (555.65,122.65) and (554.30,121.30) ..(551.60,121.30)--(503.45,121.30).. controls (500.75,121.30) and (499.40,119.95) ..(499.40,117.25)--(499.40,63.25);
\draw[mFFFFFF,line width=5.50bp] (505.20,218.50)--(493.60,218.50);
\draw[mA2461E,line width=1.80bp] (507.70,218.50)--(491.10,218.50);
\node[inner sep=0pt,anchor=base,text=mA2461E,font={\fontsize{10bp}{12bp}\selectfont}] at (508.40,226.15) {5};
\draw[mFFFFFF,line width=5.50bp] (499.40,89.30)--(499.40,77.70);
\draw[mA2461E,line width=1.80bp] (499.40,91.80)--(499.40,75.20);
\node[inner sep=0pt,anchor=base,text=mA2461E,font={\fontsize{10bp}{12bp}\selectfont}] at (489.05,70.45) {4};
\draw[mFFFFFF,line width=5.50bp] (535.40,168.95)--(535.40,157.35);
\draw[mA2461E,line width=1.80bp] (535.40,171.45)--(535.40,154.85);
\node[inner sep=0pt,anchor=base,text=mA2461E,font={\fontsize{10bp}{12bp}\selectfont}] at (546.20,166.30) {6};
\draw[mA2461E,line width=.95bp] (456.20,187.00) circle[radius=4.6bp];
\draw[mA2461E,line width=.95bp] (535.40,121.30) circle[radius=4.6bp];
\fill[mA2461E,line width=1.80bp] (426.05,131.95)--(428.10,138.55)--(424.00,138.55)--cycle;
\fill[mA2461E,line width=1.80bp] (441.05,248.20)--(434.45,250.25)--(434.45,246.15)--cycle;
\fill[mA2461E,line width=1.80bp] (434.45,53.80)--(441.05,51.75)--(441.05,55.85)--cycle;
\draw[m366448,line width=1.80bp] (331.55,217.15)--(331.55,240.55).. controls (331.55,245.65) and (329.00,248.20) ..(323.90,248.20)--(222.57,248.20).. controls (219.92,248.20) and (218.60,246.88) ..(218.60,244.23)--(218.60,238.75);
\draw[m366448,line width=1.80bp] (218.60,63.25)--(218.60,57.77).. controls (218.60,55.12) and (219.92,53.80) ..(222.57,53.80)--(323.90,53.80).. controls (329.00,53.80) and (331.55,56.35) ..(331.55,61.45)--(331.55,84.85);
\draw[m366448,line width=1.80bp] (218.60,238.75)--(218.60,191.05).. controls (218.60,188.35) and (217.25,187.00) ..(214.55,187.00)--(149.30,187.00).. controls (146.60,187.00) and (145.25,185.65) ..(145.25,182.95)--(145.25,87.55).. controls (145.25,84.85) and (146.60,83.50) ..(149.30,83.50)--(250.55,83.50).. controls (253.25,83.50) and (254.60,84.85) ..(254.60,87.55)--(254.60,141.55).. controls (254.60,144.25) and (253.25,145.60) ..(250.55,145.60)--(238.85,145.60).. controls (236.15,145.60) and (234.80,146.95) ..(234.80,149.65)--(234.80,159.10).. controls (234.80,161.80) and (236.15,163.15) ..(238.85,163.15)--(286.10,163.15).. controls (288.80,163.15) and (290.15,164.50) ..(290.15,167.20)--(290.15,214.45).. controls (290.15,217.15) and (288.80,218.50) ..(286.10,218.50)--(179.45,218.50).. controls (176.75,218.50) and (175.40,217.15) ..(175.40,214.45)--(175.40,180.70).. controls (175.40,178.00) and (176.75,176.65) ..(179.45,176.65)--(250.55,176.65).. controls (253.25,176.65) and (254.60,175.30) ..(254.60,172.60)--(254.60,154.60).. controls (254.60,151.90) and (255.95,150.55) ..(258.65,150.55)--(270.80,150.55).. controls (273.50,150.55) and (274.85,149.20) ..(274.85,146.50)--(274.85,125.35).. controls (274.85,122.65) and (273.50,121.30) ..(270.80,121.30)--(222.65,121.30).. controls (219.95,121.30) and (218.60,119.95) ..(218.60,117.25)--(218.60,63.25);
\draw[mFFFFFF,line width=5.50bp] (224.40,218.50)--(212.80,218.50);
\draw[m366448,line width=1.80bp] (226.90,218.50)--(210.30,218.50);
\node[inner sep=0pt,anchor=base,text=m366448,font={\fontsize{10bp}{12bp}\selectfont}] at (227.60,226.15) {8};
\draw[mFFFFFF,line width=5.50bp] (218.60,89.30)--(218.60,77.70);
\draw[m366448,line width=1.80bp] (218.60,91.80)--(218.60,75.20);
\node[inner sep=0pt,anchor=base,text=m366448,font={\fontsize{10bp}{12bp}\selectfont}] at (208.25,70.45) {7};
\draw[mFFFFFF,line width=5.50bp] (254.60,168.95)--(254.60,157.35);
\draw[m366448,line width=1.80bp] (254.60,171.45)--(254.60,154.85);
\node[inner sep=0pt,anchor=base,text=m366448,font={\fontsize{10bp}{12bp}\selectfont}] at (265.40,166.30) {9};
\draw[m366448,line width=.95bp] (175.40,187.00) circle[radius=4.6bp];
\draw[m366448,line width=.95bp] (254.60,121.30) circle[radius=4.6bp];
\fill[m366448,line width=1.80bp] (145.25,131.95)--(147.30,138.55)--(143.20,138.55)--cycle;
\fill[m366448,line width=1.80bp] (276.95,248.20)--(283.55,246.15)--(283.55,250.25)--cycle;
\fill[m366448,line width=1.80bp] (283.55,53.80)--(276.95,55.85)--(276.95,51.75)--cycle;
\draw[mFFFFFF,line width=5.50bp] (354.48,103.27)--(363.52,110.53);
\draw[m244C81,line width=1.80bp] (352.53,101.70)--(365.47,112.10);
\node[inner sep=0pt,anchor=base,text=m244C81,font={\fontsize{10bp}{12bp}\selectfont}] at (359.00,115.90) {1};
\draw[mFFFFFF,line width=5.50bp] (354.48,147.37)--(363.52,154.63);
\draw[m244C81,line width=1.80bp] (352.53,145.80)--(365.47,156.20);
\node[inner sep=0pt,anchor=base,text=m244C81,font={\fontsize{10bp}{12bp}\selectfont}] at (359.00,160.00) {2};
\draw[mFFFFFF,line width=5.50bp] (354.48,191.47)--(363.52,198.73);
\draw[m244C81,line width=1.80bp] (352.53,189.90)--(365.47,200.30);
\node[inner sep=0pt,anchor=base,text=m244C81,font={\fontsize{10bp}{12bp}\selectfont}] at (359.00,204.10) {3};
\node[inner sep=0pt,anchor=base,text=m366448,font={\itshape\fontsize{13bp}{15bp}\selectfont}] at (218.60,266.20) {$\operatorname{Cl}(W_2)$};
\node[inner sep=0pt,anchor=base,text=mA2461E,font={\itshape\fontsize{13bp}{15bp}\selectfont}] at (499.40,266.20) {$\operatorname{Cl}(W_2')$};
\end{tikzpicture}}
\caption{The connected sum of the marked summands in Figure~\ref{fig:virtual-summands}.}
\label{fig:insertions}
\end{figure}

We form the connected sum $K$ by cutting the diagrams at the marked
points $p_1$ and $p_2$ and joining the resulting arcs so that the orientations agree, see Figure~\ref{fig:insertions}. 
The point $p_1$ lies between $-3$ and $+1$, and $p_2$ between $+7$ and $-8$ so the Gauss word of the sum is
the concatenation 
\[
W_1W_2=(+1,-2,+3,\; W_2',\; -1,+2,-3,\; W_2), 
\]
where $W_2'=(-5,-4,-6,+5,+6,+4)$ is $W_2$ with the labels $7,8,9$ renamed $4,5,6$. 
Thus $W_1W_2$ is obtained from the word $T=(+1,-2,+3,-1,+2,-3)$ representing the right-handed trefoil by inserting $W_2'$ after $+3$
and $W_2$ after $-3$. 

\begin{figure}[h!]
\centering
\begin{tikzpicture}
\begin{scope}[xshift=-2.6cm]\gauss{1,-2,3,-5,-4,-6,5,6,4,-1,2,-3}{12}{1,...,6}
\node at (0,-2) {\(K_1\)};\end{scope}
\begin{scope}[xshift=2.6cm]\gauss{-8,-7,-9,8,9,7}{6}{7,8,9}
\node at (0,-2) {\(K_2\)};\end{scope}
\end{tikzpicture}
\caption{Marked summands. Read clockwise from the tick at the top. Each chord points from the overpass to the underpass. Intersections between chords carry no crossing data.}
\label{fig:Gauss}
\end{figure}
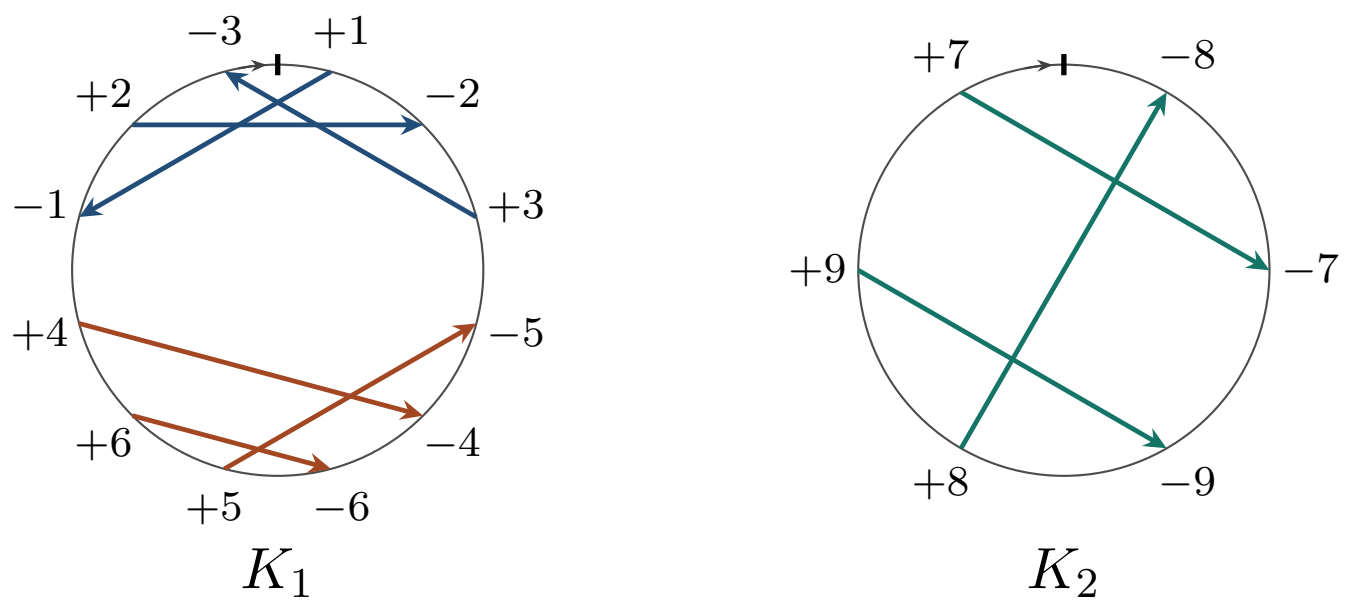

The knot $K_2$ is nonclassical and hence nontrivial.
Recall that two chords in a Gauss diagram are \emph{linked} if their endpoints
alternate around the circle. A crossing is called \emph{Gaussian odd}
if its chord is linked with an odd number of other chords. In the Gauss
diagram of $K_2$ (Figure~\ref{fig:Gauss}), each of the chords $7$ and $9$ is linked only with
chord $8$. Thus crossings $7$ and $9$ are odd, while crossing $8$ is
even. Kauffman's \emph{odd writhe} $J$ is the sum of the signs
of the Gaussian odd crossings.  In our
diagram, crossings $7$ and $9$ both have sign $+1$, so
$J(K_2)=2$.  By Kauffman's
theorem~\cite[Theorem~2(4)]{KauffmanOddWrithe}, a virtual knot
with nonzero odd writhe is nonclassical. Hence $K_2$
is nonclassical.

\section{The supporting genera}\label{sec:certificate}

In this section we show that the diagrams with Gauss words $W_1$, $W_2$, and $W_1W_2$
have minimal supporting genera $2$, $2$, and $4$, respectively. We also
determine the homology classes of the two halves at each crossing.
We use Turaev's based matrices~\cite{T}, and first recall their definition.

\subsection{Turaev's based matrices}
Let $D$ be an oriented virtual knot diagram whose classical crossings are all positive, as in our examples. Forgetting the overpasses and underpasses turns $D$ into a \emph{virtual string} $\alpha$ in the sense of Turaev, with one arrow $e_i$ in the Gauss diagram for each crossing $i$~\cite[Sec. 2.1]{T}. Turaev directs the arrow at a crossing from its tail to its head, where the
tangent vectors at the tail and at the head, both pointing along the
orientation of the knot, form a positive basis of the tangent plane in
that order~\cite[Sec.~2.2]{T}. A crossing is positive when the overpassing
tangent followed by the underpassing tangent is a positive basis~\cite[Sec.~10.2]{T}. Hence, at a positive crossing $i$ the tail of $e_i$ is the overpass $+i$ and the head is the
underpass $-i$, as in Figure~\ref{fig:Gauss}.

The canonical surface $\Sigma$ of $\alpha$ is the compact surface with boundary obtained by thickening the ribbon graph of $D$. Capping its boundary circles with disks gives the closed surface $F_D$ of Section~\ref{sec:manturov}. Let
$s\in H_1(\Sigma; \mathbb{Z})$ be the homology class of the whole curve. At crossing $i$, follow the
curve from the overpass $+i$ to the underpass $-i$. Since the two
endpoints represent the same point of $\Sigma$, this path is a loop. Denote the integral homology class of this loop by $[e_i]$. The classes $s$ and $[e_i]$, one for each crossing, form a basis of $H_1(\Sigma; \mathbb{Z})$. To see this, start with the parameter
circle and add an edge joining the two occurrences of each crossing.
Each added edge gives one new basis cycle, consisting of that edge and
the arc of the circle from $+i$ to $-i$. Now contract all the added
edges. They are pairwise disjoint, so contracting them preserves
homology. The resulting graph is the ribbon graph of $D$, and the new
basis cycles become the loops representing the $e_i$. Since $\Sigma$
retracts onto the ribbon graph, the assertion follows. 

The based matrix $T(\alpha)$ is the matrix of the intersection form $B_\alpha$ of $\Sigma$ in this basis, with $s$ first~\cite[Sec. 6.2]{T}. Since $B_\alpha$ is skew-symmetric, $T(\alpha)$ is skew-symmetric with zero
diagonal. For $e_i=(a,b)$ and $e_j=(c, d)$ with $i\ne j$, the arrow $e_j$ links $e_i$ \emph{positively} if $c\in ab$ and $d\in ba$, and $e_j$ links $e_i$ \emph{negatively} if $c\in ba$ and $d\in ab$. The entries of $T(\alpha)$ are 
\begin{equation}\label{eq:1}
B_\alpha([e_i], s)=n(e_i), \qquad B_\alpha([e_i], [e_j])=ab\cdot cd+\varepsilon,
\end{equation}
\cite[Formula (3.5.1) and Lemma 4.2.1]{T}. Here $n(e_i)$ is the number of arrows linking $e_i$ positively minus the number of arrows linking it negatively. The number $ab\cdot cd$ counts arrows with tail in the interior of the arc $ab$ and head in the interior of $cd$, minus the number of those with tail in $cd$ and head in $ab$. The correction term $\varepsilon$ is $0$ if $e_i$ and $e_j$ are unlinked, $\varepsilon=1$ if $e_j$ links $e_i$ positively, and $\varepsilon=-1$ if $e_j$ links $e_i$ negatively. 

Turaev reduces based matrices as follows \cite[Sec. 6.1]{T}. Write $r_s$ for the row corresponding to $s$  and $r_i$ for the row corresponding to $e_i$. A row $r_i$ is \emph{annihilating} if it is zero and \emph{core} if it equals $r_s$. Two distinct rows $r_i$ and $r_j$, neither of which is the base row $r_s$, are  \emph{complementary} if $r_i+r_j=r_s$. Deleting an annihilating or core row, or a complementary pair, together with the corresponding columns, is an elementary reduction. A based matrix admitting no elementary reduction is \emph{primitive}. Every based matrix reduces to a primitive one, which is unique up to relabeling \cite[Lemma 6.1.1]{T}. If two virtual strings are homotopic, their based matrices are related by elementary reductions, inverses of elementary reductions and relabeling \cite[Lemma 6.2.1]{T}. Section 6.1. of \cite{T} works with based matrices over an arbitrary abelian group, and we use it over $\mathbb{F}_2$. The reduction modulo two of an elementary reduction over $\mathbb{Z}$ is an elementary reduction over $\mathbb{F}_2$, so the above mentioned Lemma 6.2.1 of \cite{T} also holds for based matrices reduced modulo two. Indeed, if $D'$ is another diagram of the same virtual knot, the underlying
virtual strings are homotopic, so by \cite[Lemma~6.2.1]{T} the integral
based matrices of $D$ and $D'$ are joined by a chain of elementary reductions, their inverses and relabeling. Reducing every matrix
in this chain modulo two gives such a chain over $\mathbb F_2$.  

\subsection{Reduction modulo two}\label{s:4.2}

From now on homology is taken with coefficients in $\mathbb{F}_2$.
We use $s$ and $h_i$ for the images of $s$ and $[e_i]$ in $H_1(F_D)$. Since capping the boundary circles of $\Sigma$ amounts to attaching $2$-cells, the classes $s$ and $h_i$ span $H_1(F_D)$.  Let $b$ be the mod-two intersection form of $F_D$, and let $M_D$ be $T(\alpha)$ reduced modulo two, i.e., the matrix of intersection numbers of $s$ and
the $h_i$, with $s$ indexing the first row and column. Write $r_s$ and $r_i$ for the rows of $M_D$ corresponding to $s$ and $h_i$, respectively. Since the classes $s$ and $h_i$ span $H_1(F_D)$ and the intersection form on $F_D$ is
nondegenerate, a linear combination of the classes vanishes precisely
when the same linear combination of their rows vanishes. In particular,
$\operatorname{rank} M_D=\dim H_1(F_D)=2\operatorname{genus}(F_D)$, and crossing $i$ is even precisely when $r_i=0$ or $r_i=r_s$. 

Let $I_i$ be the
open interval from $+i$ to $-i$, read in the direction of the cyclic
word. Write $L_{ij}=1\in \mathbb{F}_2$
 when $I_i$ contains exactly one occurrence of $j$ and $L_{ij}=0\in \mathbb{F}_2$ otherwise. We say that a chord \emph{joins} $I_i$ to $I_j$ if it has one endpoint in $I_i$ and the other in $I_j$, but not both endpoints in $I_i\cap I_j$. 
 Reducing \eqref{eq:1} modulo $2$, we deduce that  
 $b(h_i, s)$ is the number of elements in $I_i$ modulo two (every chord linked with $e_i$ has one endpoint in $I_i$, and every other chord has zero or two), and $b(h_i, h_j)$ is $L_{ij}$ plus the number of chords joining $I_i$ to $I_j$ modulo $2$. 

 In particular, given $k\ne i, j$, if each of $I_i$ and $I_j$ contains either both endpoints of chord $k$, or neither of them, then chord $k$ does not affect $b(h_i, s)$ or $L_{ij}$, and it does not join $I_i$ to $I_j$.

\subsection{The three diagrams}

For $W_2=(-8,-7,-9,+8,+9,+7)$, the matrix $M_D$ is $4\times 4$ with rows and columns indexed by $s, h_7, h_8$ and $h_9$. The intervals $I_7,I_8,I_9$ are
$(-8)$, $(+9,+7)$ and $(+7,-8,-7)$. Their lengths modulo $2$,  together with $b(s,s)=0$, give $r_s=(0,1,0,1)$. 
Next, $L_{78}=L_{89}=1$, and $L_{79}=0$. No chord joins $I_7$ to $I_8$ or $I_9$, and only chord $7$ joins $I_8$ to $I_9$. Hence 
\[
 b(h_7,h_8)=1+0, \quad b(h_7,h_9)=0+0, \quad \textrm{and} \quad b(h_8,h_9)=1+1, 
\]
where in each sum the first summand is $L_{ij}$. This gives the matrix $B$ displayed below. 

The word $W_2'$ is obtained from $W_2$ by renaming the labels $7,8,9$ as $4, 5, 6$. Thus the matrices of $W_2$ and $W_2'$ are the same. In the right-handed trefoil word $T=(+1, -2, +3, -1, +2, -3)$ we have $I_1=(-2,+3)$,
  $I_2=(-3,+1)$ and $I_3=(-1,+2)$. Each interval is of length $2$, so $b(h_t,s)=0$. For $t\neq u$, the chords $t$
and $u$ are linked, so $L_{tu}=1$. The third trefoil chord is the
only chord joining $I_t$ to $I_u$. Hence $b(h_t,h_u)=1+1=0$, and therefore the
matrix associated with $T$ is $0$.

Recall that $W_1W_2=(+1,-2,+3,\; W_2',\; -1,+2,-3,\; W_2)$ consists of the trefoil word $T$ and two \emph{blocks} $W_2'$ and $W_2$. Each block lies between two
consecutive letters of $T$, so each trefoil interval $I$ contains either the whole block or none of it. The block $W_2'$ lies
in $I_1$ and $I_3$, and the block $W_2$ lies in
$I_2$. Within a block, every underpass precedes the corresponding
overpass, so the interval of a block chord contains every letter outside its block. By the observation at the end of Section~\ref{s:4.2}, the block chords do not affect 
the $4\times 4$ submatrix  with rows and columns indexed by $s, h_1, h_2, h_3$, so it equals the matrix of $T$. Similarly, the principal $4\times4$ submatrix with rows and columns
  indexed by $s$ and the three half classes of a block equals $B$, since every chord outside the block has both endpoints in each interval of the block. 

 For a trefoil crossing $t$ and a block crossing $k$, the chords are not linked, so $L_{tk}=0$.  A chord
of the other block has both endpoints in $I_k$, and both or neither in
$I_t$, so it does not join
$I_t$ to $I_k$. The two trefoil letters in $I_t$ belong to chords joining $I_t$ to $I_k$, so they contribute $1+1=0$. If 
  $I_t$ does not contain the block of $k$, then no other chord joins $I_t$ ot $I_k$ and $b(h_t,h_k)=0$. If it does, then both endpoints of every chord of this block lie in
$I_t$, so such a chord joins $I_t$ to $I_k$ precisely when exactly one
of its endpoints lies in $I_k$. Since every chord outside the block has
both endpoints in $I_k$, the number of these chords is $b(h_k,s)$
modulo two. Hence $b(h_t,h_k)=b(h_k,s)$.

  Finally, consider $W_1W_2$. Let $4\le i\le 6$ and $7\le j\le 9$. Then
  $I_i$ contains $W_2$ and $I_j$ contains $W_2'$ and both contain all of $T$. Hence $L_{ij}=0$, the trefoil chords do not join $I_i$ to $I_j$, a chord of $W_2'$ joins them precisely when it has exactly one endpoint in $I_i$, and a chord of $W_2$ joins them precisely when it has exactly one endpoint in $I_j$. Thus 
  $b(h_i,h_j)=b(h_i,s)+b(h_j,s)$. 

 For $W_1=(+1,-2,+3,\; W_2',\; -1,+2,-3)$, it follows that $r_1=r_3=r_s$ and $r_2=0$.
Deleting rows and columns $1,2,3$ leaves the matrix $B$ below, in the
order $s,h_4,h_5,h_6$. For $W_1W_2$, we obtain $r_1=r_3$ and
$r_1+r_2=r_s$. Thus rows $1, 3$ are core and row $2$ is annihilating for $W_1$, and rows $1, 2$ are complementary for $W_1W_2$.
Deleting rows and columns $1,2$ leaves the matrix $C$
in the order $s,h_3,h_4,h_5,h_6,h_7,h_8,h_9$. The two matrices are
\[
 B=\begin{pmatrix}
 0&1&0&1\\
 1&0&1&0\\
 0&1&0&0\\
 1&0&0&0
 \end{pmatrix},
 \qquad
 C=\begin{pmatrix}
 0&0&1&0&1&1&0&1\\
 0&0&1&0&1&0&0&0\\
 1&1&0&1&0&0&1&0\\
 0&0&1&0&0&1&0&1\\
 1&1&0&0&0&0&1&0\\
 1&0&0&1&0&0&1&0\\
 0&0&1&0&1&1&0&0\\
 1&0&0&1&0&0&0&0
 \end{pmatrix}.
\]
Both matrices have determinant $1$ in $\mathbb F_2$. Thus the
homology classes corresponding to the remaining rows are linearly independent. The row relations express each
omitted class in terms of the remaining ones, so the remaining classes
form a basis of $H_1(F_D)$ in each case. Both matrices are also primitive, since any further reduction would require a linear
dependence among their rows.

  Now let $D'$ be any diagram of the same virtual knot as $D$. Reidemeister and detour moves of virtual knot diagrams induce homotopies of the underlying virtual strings, so by Lemmas 6.1.1 and 6.2.1 of \cite{T}, applied over $\mathbb{F}_2$, $M_{D'}$ reduces to the same primitive matrix as $M_D$ up to relabeling. Deleting rows and columns does not increase rank, so
$2\operatorname{genus}(F_{D'})=\operatorname{rank}M_{D'}$ is at least
$4$, $4$ and $8$, respectively. This is the mod-two version of Turaev's genus estimate \cite[Sec. 6.3]{T}. The primitive ranks $4, 4$ and $8$ therefore give lower bounds $2, 2$ and $4$ for the supporting genera. The table below shows that our diagrams attain these bounds.  

\begin{center}
\begin{tabular}{@{}lll@{}}
\toprule
Word & Basis of $H_1(F_D;\mathbb F_2)$ & Remaining half classes\\
\midrule
$W_1$   & $s,h_4,h_5,h_6$
        & $h_1=h_3=s,\quad h_2=0$\\
$W_2$   & $s,h_7,h_8,h_9$
        & ---\\
$W_1W_2$ & $s,h_3,h_4,h_5,h_6,h_7,h_8,h_9$
         & $h_1=h_3,\quad h_2=s+h_3$\\
\bottomrule
\end{tabular}
\end{center}

\section{Computing the projection}

Recall that $K_1$ and $K_2$ are the knots represented by $W_1$ and $W_2$,
and $K$ is their specified connected sum, represented by $W_1W_2$.
By Section~\ref{sec:certificate}, all three diagrams have minimal
supporting genus, so we may therefore apply a parity step to each diagram. Recall that crossing $i$ remains classical precisely
when $h_i$ has zero image in $H_1(F_D;\mathbb F_2)/\langle s\rangle$,
that is, when $h_i=0$ or $h_i=s$. Every other classical crossing becomes
virtual, and its two occurrences are deleted from the Gauss
word~\cite[Theorem~7 and Section~4.1]{M}.

For $K_1$, the table in Section~\ref{sec:certificate} gives
$h_1=h_3=s$ and $h_2=0$. Thus crossings $1,2,3$ remain classical.
Since $s,h_4,h_5,h_6$ form a basis, crossings $4,5,6$ become
virtual. The remaining word is $(+1,-2,+3,-1,+2,-3)$. With the
given positive crossings, this is the Gauss word of the right-handed
trefoil.
Hence $g(K_1)=3_1$. For $K_2$, the classes $s,h_7,h_8,h_9$ form a basis, so all three
crossings become virtual. A diagram with only virtual crossings
reduces to a circle by detour moves. Thus $g(K_2)$ is the unknot $U$.
For $K$, we have $h_1=h_3$ and $h_2=s+h_3$. Since $s$ has zero image in the quotient, the three
half classes $h_1, h_2$ and $h_3$ have the same image in the parity group. This image
is nonzero, since $s$ and $h_3$ belong to a basis. The same basis
contains $h_4,\ldots,h_9$, so each of these classes also has nonzero
image in the quotient. Thus all nine crossings become virtual and
$g(K)=U$.

In each case the resulting knot is classical after one parity step.
For the specified connected sum we therefore have
\[
 g(K_1\mathbin{\#}K_2)=U,
 \qquad
 g(K_1)\mathbin{\#}g(K_2)=3_1\mathbin{\#}U=3_1.
\]
This gives a counterexample to Manturov's conjecture and proves Theorem~\ref{th:main}.\qed

\begin{rem}
In our example $g(K_2)$ and $g(K)$ are trivial. This is not essential. Below we explain how to construct counterexamples where all the parts/summons are nontrivial.

Let $D$ be a virtual knot diagram and let $Q$ be a classical knot.
Choose a small disk $\Delta$ in the plane, away from any marked points of
$D$, that meets $D$ in a single arc containing no crossings. Cut a
diagram of $Q$ at a point on the boundary of its unbounded region. This
turns $Q$ into a tangle in a disk with two endpoints on the boundary.
Replacing the arc in $\Delta$ by this tangle, with matching
orientations, gives a new diagram, which we call the result of
\emph{inserting $Q$ into $D$}.

Let $Q_1$ and $Q_2$ be nontrivial classical knots. Insert a diagram of
$Q_1$ into the diagram of $K_1$ and a diagram of $Q_2$ into the diagram
of $K_2$, and form the connected sum at $p_1$ and $p_2$ as before. 
The resulting diagram of the sum is the diagram of $K$ with $Q_1$ and $Q_2$ inserted. We denote the knots represented by the three new diagrams by $K_1\# Q_1$, $K_2\# Q_2$ and $K\# Q_1\#Q_2$. The last one is a specified connected sum of the first two. 
Let $D$ be any of the three diagrams of $K_1, K_2$ and $K$, and let $D'$ be the corresponding new diagram. On $F_D$ choose a small disk around each arc where an insertion is made. Making the same insertions in these disks gives a drawing of $D'$ on $F_D$. Set $F_{D'}=F_D$. We compute $M_{D'}$ using the mod two intersection pairing on $F_{D'}$. 

The inserted arc and the arc it replaces have the same endpoints and are
homotopic in the disk relative to their endpoints. Thus the class $s$
and the classes $h_i$ of all old crossings are unchanged. At each new
crossing, one of the two halves lies in the disk and hence has class
$0$, and the other half has class $s$. Hence every new crossing $q$ is even. If
$h_q=0$, then $b(h_q,x)=0$ for every $x$, so $r_q=0$. If $h_q=s$, then
$b(h_q,x)=b(s,x)$ for every $x$, so $r_q=r_s$. Thus the row of each new crossing is either $0$ or $r_s$. Deleting the new rows and columns
therefore recovers the original matrix $M_D$ by elementary reductions,
so the primitive matrices, and hence the lower bounds $2$, $2$ and $4$
for the supporting genera, are unchanged. The new diagrams are drawn on surfaces of genera $2, 2, 4$ respectively, so they attain these lower bounds. Consequently
\[
g(K_1\#Q_1)\#g(K_2\#Q_2)=3_1\#Q_1\#Q_2\neq Q_1\#Q_2=g(K\#Q_1\#Q_2)
\]
by uniqueness of prime decomposition. For example, let $Q_1=3_1$ and
let $Q_2$ be the connected sum of $n\ge1$ right-handed trefoils. The
projection of the specified sum then has $n+1$ trefoil factors, while
the connected sum of the projections has $n+2$. Distinct values of $n$
give distinct counterexamples.
\end{rem}


\end{document}